\documentclass[11pt]{article}
\usepackage{fullpage}

\usepackage{cite}
\usepackage{rotating}
\usepackage[cmex10]{amsmath}

\usepackage{algorithmic}
\usepackage{algorithm}
\usepackage{url}
\makeatletter \g@addto@macro\UrlSpecials{\do\!{\newline}}
\usepackage[hidelinks]{hyperref}
\usepackage{breakurl}
\usepackage{epsfig,epsf,psfrag,amssymb,amsfonts,latexsym,graphicx,mathrsfs,float,multirow,tabularx,diagbox}
\usepackage[dvips]{color}
\usepackage{array}
\usepackage{subfigure}

\newtheorem{theorem}{\textbf{Theorem}}

\usepackage{lipsum}

\usepackage{ifthen}
\newboolean{showcomments}
\setboolean{showcomments}{true}
\newcommand{\ye}[1]{\ifthenelse{\boolean{showcomments}}
{\textcolor{black}{#1}}{}}

\allowdisplaybreaks

\begin{document}
%
\title{Generalized Coordinated Transaction Scheduling: \\ A Market Approach to Seamless Interfaces}
%
%
%

\author{Ye~Guo,
        Yuting~Ji,
        and~Lang~Tong,
\thanks{
Ye Guo and Lang Tong are with the School of Electrical and Computer Engineering, Cornell University, Ithaca, New York, USA.(Email: \{yg299,lt35\}@cornell.edu). Yuting Ji is with the Department of Civil and Environmental Engineering, Stanford University, Stanford, California, USA.(Email: yutingji@stanford.edu).}
}

\maketitle

\begin{abstract}
A generalization of the coordinated transaction scheduling (CTS)---the state-of-the-art interchange scheduling---is proposed. Referred to as generalized coordinated transaction scheduling (GCTS), the proposed approach addresses major seams issues of CTS: the ad hoc use of proxy buses, the presence of loop flow as a result of proxy bus approximation,  and difficulties in dealing with multiple interfaces.  By allowing market participants to submit bids across market boundaries, GCTS also generalizes the joint economic dispatch that achieves seamless interchange without market participants. It is shown that GCTS asymptotically achieves seamless interface under certain conditions. GCTS is also shown to be revenue adequate in that each regional market has a non-negative net revenue that is equal to its congestion rent.  Numerical examples are presented to illustrate the quantitative improvement of the proposed approach.
\end{abstract}


%

\section{Introduction}
\subsection{Motivation}
Much of the power grid in North America is operated by independent system operators (ISOs). Each ISO is responsible for the administration of the electricity market in its operating area. Neighboring areas are connected by tie-lines, which makes it physically possible for one ISO to import power from or export power to its neighbors. It thus makes economic sense that an ISO with high generation price imports power from its neighbors that have excess and less costly resources.

The process of setting power transfer from one regional market to another, generally referred to as the \textit{interchange scheduling}, is nontrivial. If maximizing the overall system efficiency is the objective, the power flow across different operating areas should be determined by the \textit{joint economic dispatch} (JED) that treats the entire operating region as one and minimizes the overall generation cost.  But efficiency is not the only goal that governs the operations of ISOs.

ISOs in the deregulated electricity market must operate under the principles of fair representations of stakeholders (including market participants) and the financial neutrality, \textit{i.e.}, the independence with respect to traders of the market \cite[Page 152-153]{FERC2000}. As a result, ISOs rely on market participants to set the level of power transfer across market boundaries.  These market participants aim to profit from arbitrage opportunities by submitting bids to buy from one area and offers to sell in another. It is then the responsibility of ISOs who assume no financial position in the process to clear and settle these bids in a fair and transparent fashion.

As pointed out in \cite{Oren1998}, short-term operational efficiency is not always aligned with financial neutrality; there is an inherent cost associated with any market solution. A market solution to interchange scheduling creates the so-called ``seams problem'' as defined by the additional price gap over the interface compared against the seamless operation by a single ``super ISO''. The hope is that a well-designed market solution becomes more efficient as the number of market participants increases, ultimately achieving seamless interfaces, and improves the long-term performance. To our knowledge, there is no existing market solution that provably achieves the efficiency of JED. The goal of this paper is to fill, at least partially, this gap.
\subsection{Literature Review}

The interchange scheduling is a classical problem, which goes back to the 1980s \cite{EarlyHED}.  Existing solutions can be classified into two categories. The first is based on JED and aims to achieve by neighboring ISOs the best efficiency in a distributed fashion. To this end, there is an extensive literature based on primal decomposition methods \cite{Bakirtzis&Biskas:03TPS,Zhao&LitvinovZheng:14TPS,Li&Wu&Zhang&Wang::15TPS,GuoTongetc:16TPS,GuoBoseTong:17TPS} and dual decomposition methods \cite{ConejoAguado98TPS,Binetti&Etal:14TPS,Erseghe:15TPS,Chen&Thorp&Mount:03HICCS}. These methods, unfortunately, are not compatible with the existing market structure because they remove arbitrage opportunities for external market participants, in essence, requiring ISOs to trade directly with each other.

The second category includes market-based solutions that optimize the net interchange by clearing interface bids with the coordination among ISOs \cite{White&Pike:11WP,Chatterjee&Baldic:14COR}. These techniques ensure the financial neutrality of ISOs but increase the overall generation cost. The state of the art is the \textit{coordinated transaction scheduling} (CTS), which has been recently approved by FERC for implementations in ISO-NE, NYISO, PJM, and MISO \cite{FERCCTS_NENY,FERCCTS_PJMNY,FERCCTS_PJMMISO}. A key component of these methods is the use of proxy buses as trading points for external market participants. In the clearing process, power interchange is represented by injections and withdrawals at proxy buses. Consequently, the actual power flow may differ from the schedule, causing the so-called \textit{loop flow problem} \cite{LoopFlowReport}. Furthermore, CTS is limited to setting the interchange between two neighboring areas. In practice, an ISO may need to set multiple interfaces simultaneously. An extension to interchange scheduling involving more than two areas is nontrivial, see \cite{JiTong16PESGM}. Generalizations to CTS to a stochastic setting is considered in \cite{JiZhengTong:17TPS}.

\subsection{Contribution}

This paper aims to bridge the gap between the ultimate seamless solution achieved by JED and the more practical and necessarily market-based solutions. In particular, we propose a generalization of CTS, referred to as GCTS, which achieves asymptotically seamless interfaces under certain conditions. GCTS retains the structure of bidding, clearing, and settlement of CTS, thus causing no interruption of existing market operations. A key improvement over CTS is that GCTS eliminates the proxy bus approximation (thus the associated loop flow problem) inherent in all existing market-based interchange solutions. Another advantage over existing techniques is that GCTS is shown to be revenue adequate, \textit{i.e.}, the net revenue of each ISO is non-negative and is equal to its congestion rent.

GCTS can also be viewed as a generalization of JED with two modifications. First, similar to CTS but different from JED, GCTS minimizes not only the total generation cost but also the market cost of clearing interface bids. Second, similar to JED but different from CTS, GCTS solves a distributed optimization problem with all network constraints as well as an additional constraint that uses cleared interface bids to define the boundary state. This constraint is consistent with the principle of an independent and financially neutral ISO in the sense that it is the market participants who set the interchange.

GCTS does have its own shortcoming. Because GCTS solves a distributed optimization problem as in JED, it has the computation cost similar to that of JED and is more costly than CTS.  We acknowledge but do not address this issue here except pointing out that some recent techniques \cite{Zhao&LitvinovZheng:14TPS,GuoTongetc:16TPS,GuoBoseTong:17TPS} that enjoy a finite-step convergence, which alleviate to some degree such costs.


The remainder of this paper is organized as follows. We first review the CTS approach in Section II. In Section III, we present the model of multi-area power systems, the structure of interface bids, and their clearing and settlement process. Properties of GCTS are established in Section VI. Section V provides the results of simulations. 

\section{Coordinated Transaction Scheduling}

In this section, we briefly review the state-of-the-art interchange mechanism CTS \cite{White&Pike:11WP} in a deterministic setting. For a stochastic version of CTS, see \cite{JiTong16PESGM}.

For ease of presentation, we consider throughout this paper a two-area power system illustrated in Fig.\ref{fig:cts}(a). The proposed GCTS for more than two areas is straightforward and illustrated in the Appendix. We define a \textit{boundary bus} as one to which a tie-line connecting two areas is attached. Other buses are called \textit{internal buses}.

The system in Fig.\ref{fig:cts}(a) is jointly operated by ISO 1 and ISO 2. In particular, each ISO controls the interior of its operating region defined by internal buses, and the two operators control jointly operating boundaries defined by boundary buses. The interchange problem is a two-stage process in which the neighboring ISOs jointly set the boundary state in a look-ahead scheduling, and each ISO optimizes internal states in real time subject to fixed interchange schedules.

In CTS, as shown in Fig.\ref{fig:cts}(b), a ``proxy bus'' is selected among external boundary buses\footnote{In case of internal buses of the neighboring area being placed as trading locations, we can preserve these trading buses in the equivalent network on the boundary in Fig.\ref{fig:1}(b). Thereby, similar method can be derived. For simplicity, we assume hereafter that all interface bids are at boundary buses} in each area as a trading location of market participants who submit interface bids to the coordinator\footnote{It is NYISO in the implementation between New York and New England}. Each interface bid is a pair of buying and selling bids at proxy buses. They represent market participants' interest to arbitrage in a certain direction, which changes with the anticipated price gap. These bids are used to set the interchange defined as the net power transfer (rather than power flows on tie-lines) across boundaries.

\begin{figure}[ht]
  \centering
  \includegraphics[width=0.45\textwidth]{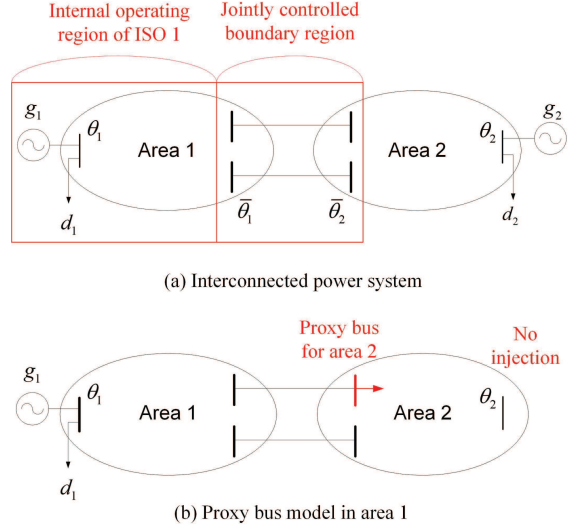}
  \caption{\small The model of the interconnected power system and the proxy-bus model in CTS}\label{fig:cts}
\end{figure}

Each interface bid has three attributes: an anticipated price difference $\Delta \pi$ at proxy buses, a maximal quantity $s_{max}$, and an import/export direction. CTS bids are cleared 15-30 minutes prior to individual real-time markets by the coordinator. The clearing process is based on the minimizing of generation cost and the payoff to market participants. To this end, the coordinator collects demand/supply curves from system operators that are used in conjunction of bids from market participants to determine the interchange quantity. The demand/supply curve from each operator is obtained by computing the expected LMP at the proxy bus for its neighboring area for each interchange level\footnote{Injections and line capacities of the neighboring area are not used}.

\begin{figure}[ht]
  \centering
  \includegraphics[width=0.5\textwidth]{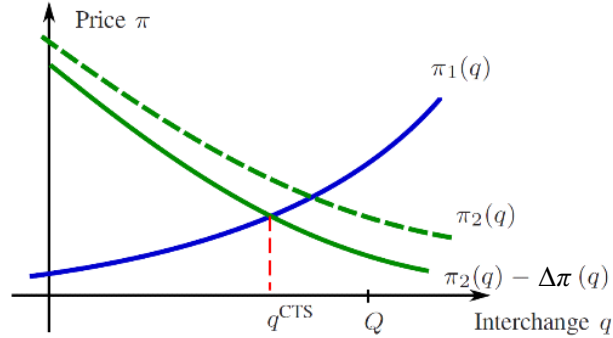}
  \caption{\small Illustration of CTS's clearing without tie-line congestion}\label{fig:ctscurve}
\end{figure}
We use the graphical representation in Fig.\ref{fig:ctscurve} from \cite{White&Pike:11WP} to illustrate the
clearing principle of CTS without tie-line congestion. Therein, curve $\pi_i (q)$
represents the incremental cost of generation for Area $i$, and
$Q$ is the interface capacity. In this example, the direction of
interchange is from Area 1 to Area 2\footnote{This is because $\pi_1(0)<\pi_2(0)$. If $\pi_1(0)>\pi_2(0)$, the direction would be opposite.}, so $\pi_1 (q)$ and $\pi_2 (q)$
serve as supply and demand curves, respectively. The third
price curve $\pi_2 (q)-\Delta \pi (q)$ is the adjusted curve of $\pi_2 (q)$ by
subtracting the aggregated interface bids $\Delta \pi (q)$. CTS interchange schedule $q^{CTS}$
is set at the intersection of $\pi_1 (q)$ and $\pi_2 (q)-\Delta \pi (q)$. All interface bids with prices
lower than $\Delta \pi (q^{CTS})$ are cleared.

Interface bids are separately settled in individual real-time
markets where the proxy bus injection is set as $q^{CTS}$.
The net interchange between the two areas will match with the scheduled $q^{CTS}$.
The real-time LMP at proxy buses $\pi_1^{RT}$ and $\pi_2^{RT}$
are used to settle cleared interface bids.
We note that there is a time latency between the clearing of
interface bids and the physical power delivery.
Such randomness may cause price deviations from the expected LMP difference at the time
of interface bid clearing. Therefore, market participants
with cleared bid offers are exposed to risks of losing money.

\begin{figure}[ht]
  \centering
  \includegraphics[width=0.5\textwidth]{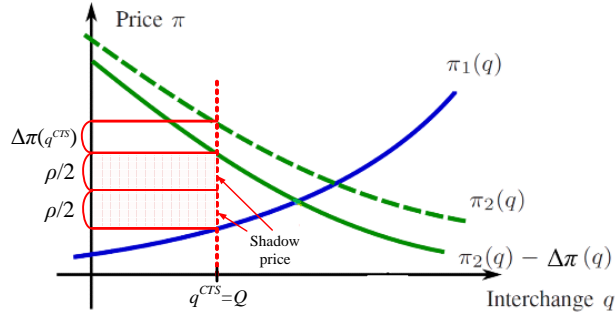}
  \caption{\small Illustration of CTS's clearing with tie-line congestion}\label{fig:ctscurveCongested}
\end{figure}

\ye{If there are tie-line congestions, \textit{i.e.}, the intersection of $\pi_1 (q)$ and $\pi_2 (q)-\Delta \pi (q)$ is greater than $Q$, then the net interchange will be scheduled at $q^{CTS}=Q$. There is $\pi_2 (q^{CTS})-\pi_1 (q^{CTS}) >\Delta \pi (q^{CTS})$. In CTS, such a price difference $\rho$ is equally partitioned into congestion prices for the two areas. Specifically, interface bids are paid at $(\pi_1^{RT}-\frac{\rho}{2})$ in Area 1 and charged at $(\pi_2^{RT}+\frac{\rho}{2})$ in Area 2, respectively. Note that $\rho$ is calculated in the look-ahead clearing process, whereas $\pi_i^{RT}$ is determined in the real-time dispatch.}

\begin{figure}[ht]
  \centering
  \includegraphics[width=0.7\textwidth]{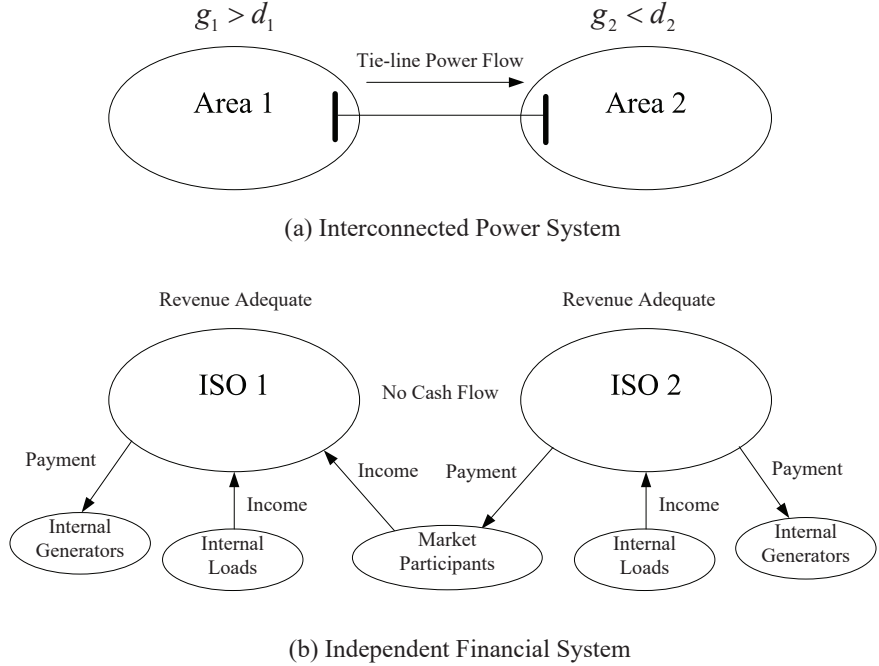}
  \caption{\small Physical power system versus the financial trading procedure.}\label{fig:illustration}
\end{figure}

We review the process of CTS and the role of interface bids via Fig. \ref{fig:illustration}. Clearing interface bids will create imbalances of local supply and demand in each area. Physically, as in Fig. \ref{fig:illustration}(a), such local imbalances naturally compel power to flow across tie-lines in a interconnected power system. Financially, as shown in Fig. \ref{fig:illustration}(b), there is no direct cash flow between the two ISOs. The payment to excess power generations in Area 1 and the revenue from excess power consumptions in Area 2 are balanced by external market participants who buy from Area 1 and sell to Area 2. Note that, cleared interface bids are financial contracts and do not physically generate or consume. They simply provide financial compensations that allow each regional market to dispatch imbalanced generations and consumptions so that power can flow across their boundaries.

Although it is reported that the CTS approach has to some extent ameliorated the seams issue, inefficient scheduling still persists \cite{NY2016Report}. In particular, modeling the net interchange as the injection to the proxy bus may be highly inaccurate when there are multiple tie-lines. In what follows, we present a generalization of CTS by removing the proxy bus approximation. 

\section{Generalized CTS}
\subsection{Network Model}

Without loss of generality, we assume that no generator or load is on the boundary bus. This assumption is made for convenience of presentation. A boundary bus that has a generator can be split into a fictitious internal bus with a generator and a boundary bus without injection. We also assume that each internal bus has one generator and one load. Let $g_i$ be the vector of power generations and $d_i$ the vector of load in Area $i$.

We adopt the DC power flow model in this paper. Specifically, nodal phase angles are state variables that are determined by active power injections. The state variables in Area $i$ are partitioned into internal phase angles $\theta_i$ and boundary phase angles $\bar{\theta}_i$. We also use $\bar{\theta}=[\bar{\theta}_{1},\bar{\theta}_{2}]^{T}$ to represent all boundary phase angles.

The DC power flow equation for the two-area power system in Fig.\ref{fig:cts}(a) is
\begin{equation}
\left[
\begin{array}{cccc}
Y_{11}&Y_{1\bar{1}}&\mbox{}&\mbox{}\\
Y_{\bar{1}1}&Y_{\bar{1}\bar{1}}&Y_{\bar{1}\bar{2}}&\mbox{}\\
\mbox{}&Y_{\bar{2}\bar{1}}&Y_{\bar{2}\bar{2}}&Y_{\bar{2}2}\\
\mbox{}&\mbox{}&Y_{2\bar{2}}&Y_{22}
\end{array}
\right]\left[
\begin{array}{c}
\theta_{1}\\
\bar{\theta}_{1}\\
\bar{\theta}_{2}\\
\theta_{2}
\end{array}
\right]=\left[
\begin{array}{c}
g_{1}-d_{1}\\
0\\
0\\
g_{2}-d_{2}
\end{array}\right],
\label{eq:ged_dclf}
\end{equation}
where $Y_{1\bar{1}}$ is the nodal admittance sub-matrix\footnote{The matrix $Y$ is composed of reciprocals of branch reactance and differs from the bus admittance matrix used in AC power flow model.} associated with the internal and boundary buses in Area 1, and $Y_{\bar{1}\bar{2}}$ the sub-matrix associated with boundary buses in areas 1 and 2.  Other terms in the coefficient matrix in (\ref{eq:ged_dclf}) are similarly defined.

\begin{figure}[ht]
  \centering
  \includegraphics[width=0.3\textwidth]{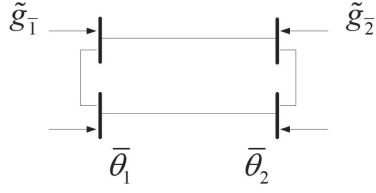}
  \caption{\small Boundary equivalent network model in GCTS}\label{fig:1}
\end{figure}

An equivalent network that captures completely the electrical properties at the boundary of the two networks can be derived as follows and is illustrated in Fig.\ref{fig:1}. Here boundary buses of Area $i$ have equivalent self-admittance $\tilde{Y}_{\bar{i}\bar{i}}$, mutual-admittance $Y_{\bar{i}\bar{j}}$, and generation $\tilde{g}_{\bar{i}}$. The power flow equation for the equivalent network for the interface is

\begin{equation}
\left[
\begin{array}{cc}
\tilde{Y}_{\bar{1}\bar{1}}&Y_{\bar{1}\bar{2}}\\
Y_{\bar{2}\bar{1}}&\tilde{Y}_{\bar{2}\bar{2}}
\end{array}
\right]\left[
\begin{array}{c}
\bar{\theta}_{1}\\
\bar{\theta}_{2}
\end{array}
\right]=\left[
\begin{array}{c}
\tilde{g}_{\bar{1}}\\
\tilde{g}_{\bar{2}}
\end{array}\right],
\label{eq:bnd_lf}
\end{equation}
where
\begin{equation}
\tilde{Y}_{\bar{i}\bar{i}}=Y_{\bar{i}\bar{i}}-Y_{\bar{i}i}Y_{ii}^{-1}Y_{i\bar{i}}, \tilde{g}_{\bar{i}}=-Y_{\bar{i}i}Y_{ii}^{-1}(g_i - d_i).\label{eq:yeq}
\end{equation}

The coefficient matrix in (\ref{eq:bnd_lf}) does not change with nodal power injections. Throughout this paper, we assume that the two-area system is on the same island, so the coefficient matrix in (\ref{eq:bnd_lf}) is full rank after removing the reference bus. The equivalent power injection $\tilde{g}_{\bar{i}}$ succinctly captures the external impact of internal power injections in Area $i$; it represents its power interchange schedule. Therefore, in what follows, those equivalent power injections are associated with interface bids from external market participants. Hereafter, we drop the word ``external'' if that does not cause any confusion.

\subsection{Definition of Interface Bids}

GCTS uses the same format of bids as CTS.  Namely, an interface bid $i$ is defined by a triple
\begin{equation}
\mathcal{B}\triangleq\{<B_{pm},B_{qn}>, \Delta \pi_i, s_{\rm max,i}\},\nonumber
\end{equation}
where
\begin{enumerate}
\item $<B_{pm},B_{qn}>$ is an ordered pair of boundary buses that specifies the bid as withdrawing at bus $m$ in Area $p$ and injecting the same amount at bus $n$ in Area $q$. They need not be directly connected by a tie-line;
\item $\Delta \pi_i$ is its price bidding on the anticipated price gap that the bid is settled in the two real-time markets\footnote{This may not be equal to the LMP difference. See Subsection III-D and Remark 2 after Theorem \ref{thm:JED} for mathematical and economical interpretations.};
\item $s_{\rm max,i}$ is its maximum quantity.
\end{enumerate}

\begin{figure}[ht]
  \centering
  \includegraphics[width=0.29\textwidth]{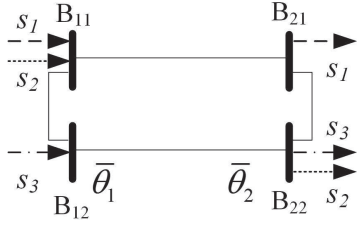}
  \caption{\small Network equivalence on the boundary. Dotted-line arrows represent three interface bids in the example below: $s_1$ injects at $B_{11}$ and withdraws at $B_{21}$; $s_2$ injects at $B_{11}$ and withdraw at $B_{22}$; $s_3$ injects at $B_{12}$ and withdraws at $B_{22}$.}\label{fig:net_eq}
\end{figure}

The only difference between CTS and GCTS is that, in stead of using a single proxy bus in each area,  GCTS allows bids to be submitted to all pairs of boundary buses across the boundary, as illustrated in Fig \ref{fig:net_eq}.

We aggregate all interface bids with an incidence matrix $M_i$ associated with boundary buses in Area $i$. Specifically, each row of $M_i$ corresponds to a boundary bus of Area $i$, and each column of which corresponds to an interface bid. The entry $M_{i}(m,k)$ is equal to one if interface bid $k$ buys power at boundary bus $B_{im}$ from Area $i$, minus one if it sells power at bus $B_{im}$ to Area $i$, and zero otherwise. For example, if there are three bids as illustrated in Fig.\ref{fig:net_eq}, matrices $M_i (i=1,2)$ are
\begin{equation}
M_1\!=\!
\left[\!
\begin{array}{ccc}
1\!&\!1\!&\!\!0\\
0\!&\!0\!&\!\!1
\end{array}\!
\right] \!\!\!\begin{array}{c}
(B_{11})\\
(B_{12})
\end{array}, M_2\!=\!
\left[\!
\begin{array}{ccc}
-1\!\!&\!0\!&\!0\\
0\!\!&\!\!-1\!\!&\!\!-1
\end{array}\!
\right]\!\!\! \begin{array}{c}
(B_{21})\\
(B_{22})
\end{array}.
\end{equation}

Consequently, let $s$ be the vector whose $i$th entry $s_i$ is the cleared quantity of bid $i$. Then $M_i s$ represents the aggregated equivalent power injection induced by cleared interface bids on boundary buses in Area $i$. By substituting the right-hand side in (\ref{eq:bnd_lf}) by $M_i s$, we have

\begin{equation}
\left[
\begin{array}{cc}
\tilde{Y}_{\bar{1}\bar{1}}&Y_{\bar{1}\bar{2}}\\
Y_{\bar{2}\bar{1}}&\tilde{Y}_{\bar{2}\bar{2}}
\end{array}
\right]\hspace{-0.1cm}\left[\!
\begin{array}{c}
\bar{\theta}_{1}\\
\bar{\theta}_{2}\\
\end{array}\!
\right]\!=\!\left[\!
\begin{array}{c}
M_1 s\\
M_2 s
\end{array}\!\right].
\label{eq:bnd_dclf}
\end{equation}

In (\ref{eq:bnd_dclf}), the interchange schedule is solely determined by the cleared interface bids from market participants. In the market clearing process of GCTS, as presented in the next subsection, Equation (\ref{eq:bnd_dclf}) will be incorporated as an equality constraint in the optimization model where the internal bids $g_i$ and interface bids $s$ are cleared together.

\subsection{Market Clearing Mechanism}
GCTS preserves the architecture of CTS; it assumes the presence of a coordinator who collects interface bids and clears them via a look-ahead dispatch, and the interface bids are settled separately in the real-time markets. GCTS removes the proxy bus approximation, and its clearing of interface bids is based on a generalization of JED. The key idea is to clear interface bids by optimizing the boundary state as follows:

\begin{align}
\min \limits_{\{g_{i},s,\bar{\theta}, \theta_i\}}\! & c(g_1,g_2,s)=\sum \limits_{i=1}^{2}c_{i}(g_{i})+ \Delta \pi^T s, \label{eq:mc_obj}\\
\textrm{subject to}&\hspace{0.1cm}\check{g}_i\leq g_i\leq \hat{g}_i, i=1,2,\label{eq:glimits}\\
&\hspace{0.1cm}0\leq s\leq s_{\textrm{max}}, \label{eq:slimits}\\
&\hspace{0.1cm}H_{i}\theta_{i}+H_{\bar{i}}\bar{\theta}_{i}\leq f_{i},i=1,2, \label{eq:mc_linecons}\\
&\hspace{0.1cm}\bar{H}_{\bar{1}}\bar{\theta}_{1}+\bar{H}_{\bar{2}}\bar{\theta}_{2}\leq \bar{f}, \label{eq:mc_tieline}\\
&\hspace{-2.0cm} \left[\!
\begin{array}{cccc}
Y_{11}&Y_{1\bar{1}}&\mbox{}&\mbox{}\\
Y_{\bar{1}1}&Y_{\bar{1}\bar{1}}&Y_{\bar{1}\bar{2}}&\mbox{}\\
\mbox{}&Y_{\bar{2}\bar{1}}&Y_{\bar{2}\bar{2}}&Y_{\bar{2}2}\\
\mbox{}&\mbox{}&Y_{2\bar{2}}&Y_{22}
\end{array}\!
\right]\!\!\left[\!
\begin{array}{c}
\theta_{1}\\
\bar{\theta}_{1}\\
\bar{\theta}_{2}\\
\theta_{2}
\end{array}\!
\right]\!=\!\left[\!
\begin{array}{c}
g_{1}-d_{1}\\
0\\
0\\
g_{2}-d_{2}
\end{array}\!\right], \label{eq:mc_pf}\\
& \hspace{-0.5cm} \left[
\begin{array}{cc}
\tilde{Y}_{\bar{1}\bar{1}}&Y_{\bar{1}\bar{2}}\\
Y_{\bar{2}\bar{1}}&\tilde{Y}_{\bar{2}\bar{2}}
\end{array}
\right]\!\left[\!
\begin{array}{c}
\bar{\theta}_{1}\\
\bar{\theta}_{2}\\
\end{array}\!
\right]\!=\!\left[\!
\begin{array}{c}
M_1 s\\
M_2 s
\end{array}\!\right],\label{eq:mc_bndpf}
\end{align}
where decision variables are the cleared internal generation bids $g_i$ with the quantity limit (\ref{eq:glimits}), the cleared interface bids $s$ with the quantity limit (\ref{eq:slimits}), and the system states ($\theta_i, \bar{\theta}$) subject to internal and tie-line power limits (\ref{eq:mc_linecons}) and (\ref{eq:mc_tieline}). Any bid $i$ with $s_i=s_{\rm max, i}$ is fully cleared, any with $s_i=0$ is rejected, and any with $0<\!s_i\!<s_{\rm max, i}$ is partially cleared at amount $s_i$.

Note that, the term $\Delta\pi^T s$ represents the market cost of clearing interface bids. Because the price difference in the real time is in general different from the look-ahead dispatch, market participants carry a certain amount of risk.  Thus the bid $\Delta\pi_i$ represents the willingness of the bidder $i$ to take that risk.  See \cite{MEAN_VAR_1995} for details of the quantification for risks.

The market clearing model of GCTS (\ref{eq:mc_obj})-(\ref{eq:mc_bndpf}) differs from JED in two aspects: (i) the market cost of clearing interface bids $\Delta\pi^T s$ in the objective function and (ii) the additional equality constraint (\ref{eq:mc_bndpf}) that determines the boundary state by clearing interface bids subject to their quantity limits. The coordinator sets the interchange by clearing the interface bids to minimize the overall cost subject to operational constraints and constraints (\ref{eq:slimits}) and (\ref{eq:mc_bndpf}) imposed by the interface bids.

The clearing problem (\ref{eq:mc_obj})-(\ref{eq:mc_bndpf}) of GCTS is a look-ahead economic dispatch where the load powers $d_i$ are predicted values. It should be solved in a hierarchical or decentralized manner. Any effective multi-area economic dispatch method can be employed. See, \textit{e.g.}, \cite{GuoTongetc:16TPS,Zhao&LitvinovZheng:14TPS} where (\ref{eq:mc_obj})-(\ref{eq:mc_bndpf}) is solved with a finite number of iterations.

\subsection{Real-time Dispatch and Settlement}

Interface bids are settled in the real-time market together with internal bids. There is no coordination required at this step. Specifically, ISO 1 solves its local economic dispatch with fixed boundary state $\bar{\theta}$:

\begin{eqnarray}
\min \limits_{\{g_{1},\theta_1\}} & c_{1} (g_{1}), & \label{eq:rt_obj} \\
\textrm{subject to}&\hspace{0.1cm}H_{1}\theta_{1}+H_{\bar{1}}\bar{\theta}_{1}\leq f_{1}, &(\eta_1^R) \label{eq:rt_line} \\
&\hspace{0.1cm}\check{g}_1\leq g_1\leq \hat{g}_1, i=1,2,& (\bar{\xi}_1^R, \underline{\xi}_1^R)\label{eq:rt_glimits}\\
&\hspace{-2.5cm}
\left[\!
\begin{array}{ccc}
Y_{11}&Y_{1\bar{1}}&\mbox{}\\
Y_{\bar{1}1}&Y_{\bar{1}\bar{1}}&Y_{\bar{1}\bar{2}}
\end{array}\!
\right]\!\!\left[\!
\begin{array}{c}
\theta_1\\
\bar{\theta}_1\\
\bar{\theta}_2
\end{array}\!
\right]\!\!=\!\!\left[\!
\begin{array}{c}
g_1-d_1^R\\
0
\end{array}\!\right],\!\!&\!\!\begin{array}{c}
(\lambda_1^R)\\
(\bar{\lambda}_1^R)
\end{array}\label{eq:rt_lf1}
\end{eqnarray}
where $d_1^{R}$ represents real-time internal loads, which may deviate from their predictions in the look-ahead dispatch (\ref{eq:mc_obj})-(\ref{eq:mc_bndpf}). The real-time internal dispatch in each area should be compliant with the pre-determined interchange schedule. To this end, boundary state $\bar{\theta}$ is fixed at the solution to Equation (\ref{eq:mc_bndpf}) with $s$ cleared interface bids solved from (\ref{eq:mc_obj})-(\ref{eq:mc_bndpf}). All multipliers are given to the right of corresponding constraints.

ISO 1 simultaneously settles internal and interface bids in the real-time market. Internal bids are settled at the LMP $\lambda_1^R$. To settle interface bids, we need to analyze the sensitivity of the local optimal cost in (\ref{eq:rt_obj}) with respect to $s$. In the real-time dispatch (\ref{eq:rt_obj})-(\ref{eq:rt_lf1}), the impact of interface bids $s$ is imposed via the fixed boundary state variables $\bar{\theta}$. The sensitivity of local optimal cost with respect to $\bar{\theta}$ is
\begin{equation}
\nabla_{\bar{\theta}} c_1^*=
\left[\!
\begin{array}{cc}
Y_{\bar{1}1}&Y_{\bar{1}\bar{1}}\\
            &Y_{\bar{2}\bar{1}}
\end{array}\!
\right]\left[\!
\begin{array}{c}
\lambda_1^R \\
\bar{\lambda}_1^R
\end{array}\!
\right]+\left[\!
\begin{array}{c}
H_{\bar{1}}^T \eta_1^R \\
0
\end{array}\!
\right].
\end{equation}

The sensitivity of local optimal cost with respect to $s$ is

\begin{equation}
\nabla_{s} c_1^* \!=\! [\nabla_{s} \bar{\theta}]^T \nabla_{\bar{\theta}} c_1^* \!=\! M^T \left[\!
\begin{array}{cc}
\tilde{Y}_{\bar{1}\bar{1}}&Y_{\bar{1}\bar{2}}\\
Y_{\bar{2}\bar{1}}&\tilde{Y}_{\bar{2}\bar{2}}
\end{array}\!
\right]^{-1}\!\!\! \nabla_{\bar{\theta}} c_1^*\triangleq\! \mu_1^R. \label{eq:rt_mu}
\end{equation}

In the absence of tie-line congestion, interface bids pay prices $\mu_1^R$ in Area 1 and $\mu_2^R$ in Area 2 (they get paid if $\mu_i<0$). In general, interface bids are not settled at LMPs. This is because the change of the objective function \eqref{eq:rt_obj} with an increment of cleared $s$ differs from that with an increment of load power.

If there are tie-lines congested, similar to CTS, we will compute congestion rents according to the look-ahead dispatch (\ref{eq:mc_obj})-(\ref{eq:mc_bndpf}) and subtract them from the payment to interface bidders. Tie-line congestion prices associated with interface bids are calculated by
\begin{equation}
\rho=M^T \tilde{S}^T \bar{\eta}, \tilde{S}= [\bar{H}_{\bar{1}} \hspace{0.1cm} \bar{H}_{\bar{2}}]\left[
\begin{array}{cc}
\tilde{Y}_{\bar{1}\bar{1}}&Y_{\bar{1}\bar{2}}\\
Y_{\bar{2}\bar{1}}&\tilde{Y}_{\bar{2}\bar{2}}
\end{array}
\right]^{-1},
\label{eq:bndcongestion}
\end{equation}
where $\bar{\eta}$ is the shadow price in (\ref{eq:mc_tieline}), and $\tilde{S}$ is the shift factor of boundary buses with respect to tie-lines in Fig. \ref{fig:net_eq}. Similar to CTS, we evenly split the congestion rent price $\rho$ into two areas\footnote{When there are more than two areas, tie-line congestions may induce positive shadow prices $\rho$ for interface bids over other interfaces. Nevertheless, the calculation of $\rho$ is the same as in \eqref{eq:bndcongestion}, and the shadow price should be evenly split by neighboring areas.}. Namely, market participants pay $\mu_i^{R}+\frac{\rho}{2}$ in Area $i, i=1,2$.

If $d_i^R=d_i$, one can prove that the real-time dispatch level and prices are consistent with the look-ahead dispatch. Note that fixing some variables at their optimal values does not change optimal values of other primal and dual variables. If the real-time dispatch (\ref{eq:rt_obj})-(\ref{eq:rt_lf1}) is infeasible, ad hoc adjustments such as relaxations of flow limits can be employed in practice. 

\section{Properties of GCTS}

\subsection{Efficiency and price convergence of GCTS}

By removing the proxy bus approximation and adopting a strict DC OPF model in (\ref{eq:mc_obj})-(\ref{eq:mc_bndpf}), we are able to establish many important properties for GCTS. Their proofs are relegated to the Appendix.

We first show that GCTS asymptotically achieves seamless interfaces when more and more bidders participate in the competition at all possible pairs of trading locations. Intuitively, for GCTS to achieve the cost of JED, two conditions are necessary in general. First, there have to be enough bidders who try to capture the arbitrage profits across the interface so that they drive $\Delta \pi \rightarrow 0$.  This follows the standard economic argument of perfect competition. Second, bids need to be diverse enough to make the matrix $M$ full row rank so that the tie-line flows of the GCTS can match those of JED. It turns out that both conditions can be satisfied simultaneously by conditions below.

\begin{theorem}
\label{thm:JED}
(Asymptotic efficiency) Consider a market with $N$ independent interface bidders. Assume that (i) both JED and GCTS are feasible and each has an unique optimum, (ii) the number of aggregated bids for each pair of source and sink buses grows unbounded with $N$, and (iii) bidding prices for all participants go to zero as $N \rightarrow \infty$, \textit{i.e.} $\lim_{N\rightarrow \infty} \Delta \pi =0$, then the scheduled tie-line power flows and generations in each area by GCTS converge to those of the JED as $N \rightarrow \infty$.

\end{theorem}

\textbf{Remark 1}: Recall that JED by a ``super ISO'' provides the lowest possible generation cost, thus achieving the overall market efficiency. In practice, the power system is artificially partitioned into multiple subareas that are operated by financially neutral ISOs, and interchange scheduling has to rely on bids from market participants. Such operational regulations will naturally create seams at interfaces. Theorem \ref{thm:JED} shows that, however, GCTS asymptotically achieves seamless interfaces under mild conditions. This indicates that GCTS leads to the price convergence between regional electricity markets.

\textbf{Remark 2}: The price convergence implies that there is no arbitrage opportunity, and that the dispatch level of GCTS is the same as that of JED. Note that due to congestions, boundary buses may have different LMPs even under the administration of a ``super-ISO''. So the price convergence is in fact for shadow prices of $s$. This also explains why interface bids should be settled at $\mu_i^R$ in \eqref{eq:rt_mu} but not LMPs.

\textbf{Remark 3}: The assumptions that bidding locations are diverse enough and that $\Delta\pi$ goes to zero as $N$ increases come from the interpretation that, as the number of bidders increases, there are always enough bids that can be cleared to satisfy the desired interchange level. Thus individually, each bidder seeks trading locations with seams and reduces its bidding price so that it will have a better chance to be cleared.

\subsection{Relation Between GCTS and CTS}
Next we establish connections between GCTS and CTS.  Specifically, we show that the two mechanisms are equivalent in a particular simple setting.

\begin{theorem}
\label{thm:CTS}
When there is a single tie-line between two areas, the clearing process of GCTS (\ref{eq:mc_obj})-(\ref{eq:mc_bndpf}) provides the same interchange as that of CTS.
\end{theorem}

\textbf{Remark:} A natural corollary of Theorem \ref{thm:CTS} is that when there is a single tie-line between two areas and real-time load is the same as the load considered in the interchange scheduling, then CTS provides the optimal interchange schedule in the sense that the posterior real-time dispatch $g_i^R$ minimizes the total cost of all internal and external market participants.

In practice, however, neither condition in these two theorems is likely to hold. In such cases, our simulations show that GCTS generally has lower overall cost than CTS and its dispatch satisfies security constraints. CTS, may on the other hand, may violate security constraints due to the loop flow problem engendered by its proxy-bus approximation. See Section V for details.

\subsection{Revenue Adequacy}

In this subsection, we establish the revenue adequacy for the real-time market (\ref{eq:rt_obj})-(\ref{eq:rt_lf1}). Recall that, in the single-area economic dispatch, each area has a non-negative net revenue, which is equal to its congestion rent. We prove in the following theorem that each area achieves its revenue adequacy in the same fashion in GCTS in an interconnected power system.

\begin{theorem}
\label{thm:Revenue}
Assume that the real-time dispatch (\ref{eq:rt_obj})-(\ref{eq:rt_lf1}) is feasible and that the settlement process follows our description in Subsection III-D, then
the net revenue of each area is non-negative and is equal to its congestion rent.
\end{theorem}

\subsection{Local Performance}

ISOs are mainly responsible of the efficiencies of their own regional markets, rather the overall efficiency. Therefore, an ISO may be reluctant to implement any interchange scheduling approach that worsens its local performance for the sake of the overall efficiency. We partly address this issue in this subsection.

In the conventional interchange scheduling before CTS, market participants split their bidding prices into $\Delta \pi=\pi_1+\pi_2$ and separately submit them to the two neighboring ISOs who clear these bids independently. Only bids cleared in both markets will be scheduled \cite{NENY:BeforeCTS}. In essence, we take the minimum of the cleared quantities. In the following theorem, we prove that GCTS achieves higher local surpluses in all areas than the conventional approach under a simple setting:

\begin{theorem}
\label{thm:Surplus}
Assume that (i) there is a single tie-line between two neighboring areas, (ii) real-time load demands are the same as their look-ahead predictions, and (iii) each market clearing problem has an unique optimum, then there is
\begin{equation}
\tilde{LS}_i \geq \hat{LS}_i,
\end{equation}
where $\tilde{LS}_i$ is the local surplus of area $i$ in its real-time market (\ref{eq:rt_obj})-(\ref{eq:rt_lf1}) with $\bar{\theta}$ determined by the optimal $\tilde{s}$ cleared in GCTS (\ref{eq:mc_obj})-(\ref{eq:mc_bndpf}). Specifically, it is defined as
\begin{equation}
\tilde{LS}_i \triangleq (D_i-(\tilde{\lambda}_i^R)^T d_i)+((\tilde{\lambda}_i^R)^T \tilde{g}_i - c_i(\tilde{g}_i))+f_1^T \tilde{\eta}_1^R,
\end{equation}
where $D_i$ is the constant utility of consumers. Variables with tildes are solved with $\tilde{s}$. The local total surplus in Area $i$ is the sum of its consumer surplus, supplier surplus, and the surplus of transmission owners. The local surplus $\hat{LS}_i$ with $\hat{s}$ the result of separate clearing is similarly defined. 
\end{theorem}

We remove this result from our journal submission because this is more about CTS. In general, when there are multiple tie-lines, Theorem \ref{thm:Surplus} may not hold for GCTS. Nevertheless, it is important to look into performances in regional markets, especially for power systems that cover multiple regions or even countries. Investigating weaker conditions for Theorem \ref{thm:Surplus} would be an interesting direction for future works.

\section{Numerical Tests}
\subsection{Two-area 44-bus system}
\begin{figure}[ht]
  \centering
  \includegraphics[width=0.6\textwidth]{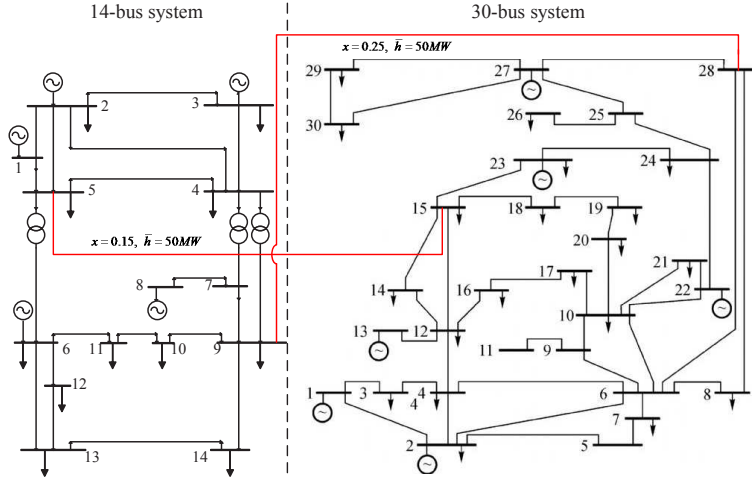}
  \caption{\small Configuration of the two-area power system}
  \label{fig:sys}
\end{figure}

GCTS was tested on a two-area system composed of the IEEE 14-bus system (Area 1) and the 30-bus system \cite{Case30} (Area 2). The system configuration and reactance and capacities of tie-lines are illustrated in Fig.\ref{fig:sys}.

Two groups of simulations were conducted. First, we aimed to illustrate the market clearing process of GCTS. Second, we compared GCTS with JED and CTS and numerically demonstrated the asymptotic convergence of GCTS to JED as in Theorem \ref{thm:JED}.

\subsubsection{Illustration of the market clearing process}

Eight interface bids were considered in the first group of simulations. Their trading locations and prices $\Delta \pi$ are listed in TABLE \ref{table:bid}. Some market participants traded on boundary buses without direct connections, such as bids 2 and 5. The maximal quantities of all bids were set as 30MW.

From default prices in Area 1, we used a weighting factor $w$ to generate scenarios with various degrees of price discrepancies. For all scenarios, cleared interface bids, tie-line power flows, marginal prices, and system costs are presented in TABLE \ref{table:result}:

\begin{table}[ht]
\centering
\caption{\small Profile of Interface Bids}
\begin{tabularx}{11cm}{cccc}
\hline
Indices&Sell to&Buy from&Price (\$/MWh)\\
\hline
1& Bus 15 (Area 2) & Bus 5 (Area 1) & 1\\
2& Bus 28 (Area 2) & Bus 5 (Area 1) & 2\\
3& Bus 5 (Area 1) & Bus 15 (Area 2) &1.5\\
4& Bus 5 (Area 1) & Bus 28 (Area 2) &0.5\\
5& Bus 15 (Area 2)& Bus 9 (Area 1) &1.0\\
6& Bus 28 (Area 2)& Bus 9 (Area 1) &2.0\\
7& Bus 9 (Area 1) &Bus 15 (Area 2)&1.5\\
8& Bus 9 (Area 1) &Bus 28 (Area 2)&0.5\\
\hline
\end{tabularx}\label{table:bid}
\end{table}

\begin{table}[ht]
\centering
\caption{\small Results of Interface Bid Clearing}
\begin{tabularx}{11cm}{cccccc}
\hline
\multicolumn{2}{c}{Weighting factor}&0.1&0.15&0.2&1.0\\
\hline
\multirow{8}{*}{\shortstack{Cleared\\quantities\\of\\interface\\bids\\(MW)}}&1& 30    & 30     & 5.56     & 5.56 \\
&2& 0    & 0     & 0     & 0  \\
&3& 0    & 0     & 0     & 0\\
&4& 0   & 0 & 30     & 30\\
&5& 0   & 0    & 0    & 0\\
&6& 0 & 0    & 0    & 0\\
&7& 10.40 & 30 & 30 & 30\\
&8& 30    & 30     & 30     & 30\\
\hline
Tie-line & bus 15 to 5 & -8.66 & 5.53  & 38.60 & 38.60 \\
  flow (MW)   & bus 28 to 9 & 19.05  & 41.31 & 45.84 & 45.84 \\
\hline
\multirow{4}{*}{\shortstack{Marginal \\prices\\(\$/MWh)}}& bus 5 & -0.10  & 0.02 & 0.82 & 15.62 \\
& bus 9 & 2.90 & 4.34 & 6.49 & 45.96\\
& bus 15  &  1.40 & 1.04 & 1.82 & 16.62\\
& bus 28  & 0.0 & 0.0 & 0.0 & 0.0\\
\hline
\multirow{3}{*}{\shortstack{Market \\costs\\(\$/h)}} & Internal & 923.2 & 1148.2 & 1371.0 & 4525.4 \\
& Interface & 58.1   & 90 & 80.56 & 80.56  \\
& Total & 983.0      & 1238.2 & 1451.6 & 4605.9\\
\hline
\end{tabularx}\label{table:result}
\end{table}

The second block (second to ninth rows) in TABLE \ref{table:result} lists cleared amounts for interface bids. With the increase of $w$, bids delivering power from Area 2 to Area 1 were cleared at greater quantities, see the fifth and eighth rows, while those delivering power in the opposite direction were cleared at smaller quantities, see the second row.

The third block includes results on tie-line power flows. They were determined by the boundary power flow equation (\ref{eq:bnd_dclf}) and cleared amounts of bids in the second block.  When $w=0.1$, tie-line power flows were in both directions. When $w$ was increased, which signified greater price discrepancies, tie-line power flows became unidirectional from the low-price area to the high-price area.

The fourth block are marginal prices for all boundary buses in the market clearing process, \textit{i.e.}, multipliers associated with boundary equality constraint (\ref{eq:mc_bndpf}). All bids whose prices were lower than marginal price gaps between their trading points were totally cleared, see the second, twelfth, and fourteenth rows when $w=0.1$ and the second row in TABLE \ref{table:bid} as an example. All bids whose prices were higher were rejected, see the third, twelfth, and fifteenth rows when $w=0.1$ and the third row in TABLE \ref{table:bid} as an example. For partially cleared interface bids, marginal price gaps between their trading points were equal to their bidding prices, see the eighth, thirteenth, and fourteenth rows when $w=0.1$ and the eighth row in TABLE \ref{table:bid} as an example.

The last block are generation costs, costs of market participants, and total costs per hour in the proposed approach. GCTS considered the total market cost of internal and interface bidders.

\subsubsection{Comparison with Existing Benchmarks}

In the second group of simulations, we compared the proposed method with existing approaches on tie-line scheduling. Specifically, the following methods were compared:

i) JED that minimized the total generation cost;

ii) CTS wherein proxy buses were selected as bus 5 in Area 1 and bus 15 in Area 2;

iii) The proposed mechanism of GCTS.

Default generation prices were considered in this test. We used similar bids to those in TABLE \ref{table:bid} for GCTS but their quantity limits and prices were uniformly set as $s_{\textrm max}=100MW$ and $\Delta \pi=\$0.1/MWh$. In CTS, all bids were placed at proxy buses with the same quantity limits and prices.

We compared market costs in the look-ahead interchange scheduling as well as those in the real-time local dispatch. For the latter, we generated 100 normally distributed realizations of real-time load consumptions whose mean values were their look-ahead predictions (default values in the system data) and standard deviations were 5\% of their mean values. Comparisons on net interchange quantities, look-ahead generation costs and total costs, and real-time average total costs for all samples are recorded in TABLE \ref{table:compare1}:

\begin{table}[ht]
\centering
\caption{\small Comparison of JED, CTS, and GCTS for the Two-area Test}
\begin{tabularx}{11cm}{cccc}
\hline
                    & JED     & CTS    & GCTS \\
\hline
Net interchange amounts (MW) & 87.0    & 80.3   & 87.0 \\
Look-ahead generation costs (\$/h)          & 4039.8  & 4109.9 & 4039.8\\
Look-ahead total costs (\$/h)         & --      & 4118.0 & 4048.5\\
Average real-time total costs (\$/h)          & 4096.2  & 4139.8 & 4115.7\\
\hline
\end{tabularx}\label{table:compare1}
\end{table}

From TABLE \ref{table:compare1} we observed that GCTS achieved lower look-ahead and average real-time costs than CTS. Specifically, GCTS had lower real-time costs in 88 out of the 100 samples. In addition, CTS suffered from the loop-flow problem in that branch power flows solved with the global power flow equation and real-time dispatch levels in both areas deviated from internal real-time schedules. In this test, average discrepancies on tie-line power flows were 18.25\% for Area 1 and 16.38\% for Area 2, respectively. As a result, CTS caused unpredicted overflows for transmission lines in all of the 100 scenarios, with 2.72 overflowed transmission lines in each scenario on average and the average ratio of overflows as 11.27\%. In GCTS, however, such problems did not exist because it is based on the exact DC power flow model. Another takeaway of TABLE \ref{table:compare1} is that, with sufficient bids and relatively low prices ($\Delta \pi=\$0.1/MWh$), the interchange scheduled by GCTS was the same as that in JED in this test.

\begin{figure}[ht]
  \centering
  \includegraphics[width=0.6\textwidth]{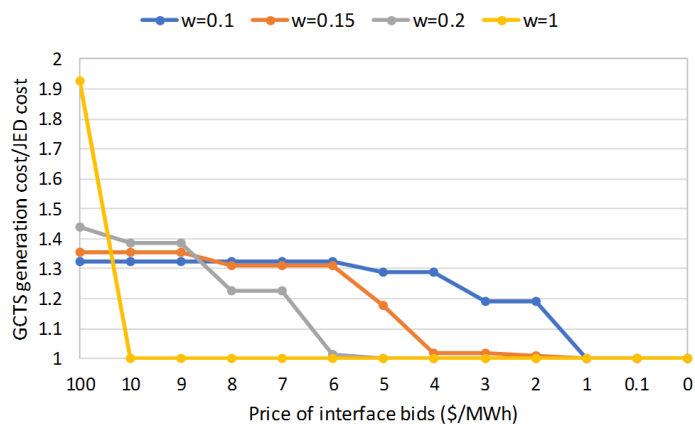}
  \caption{\small Price convergence of GCTS with different bidding prices}
  \label{fig:priceconvergence}
\end{figure}

We illustrate the price convergence of GCTS with different values of $w$ in Fig. \ref{fig:priceconvergence} by adjusting the uniform bidding price $\Delta \pi$. No bid was cleared when the bidding price $\Delta \pi=\$100/MWh$. When $\Delta \pi$ decreased to small enough values (\$0.1/MWh in this test), generation costs of GCTS in all scenarios were equal to those of JED. In general, the more significant the price discrepancy was, the faster the price converged. This is consistent with our intuition that market participants could be cleared at higher prices when there is more room for arbitrations.

Note that such price convergence did not happen in CTS. For the test in TABLE \ref{table:compare1}, for example, if we set the bidding price of CTS as zero, the total generation cost would be \$4109.7 per hour, which was higher than that of JED.

\subsection{Three area 189-bus system test}
\begin{figure}[ht]
  \centering
  \includegraphics[width=0.6\textwidth]{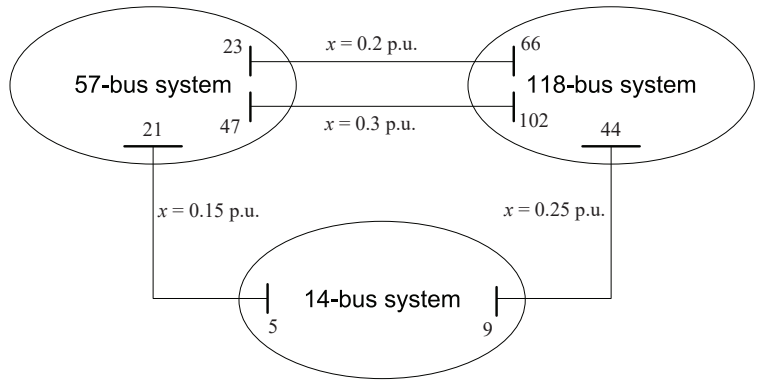}
  \caption{\small Configuration of the three-area power system}
  \label{fig:sys2}
\end{figure}

The proposed method was also tested on a three-area system as shown in Fig. \ref{fig:sys2}. The system was composed of IEEE {14,} 57, and 118-bus systems. Power flow limits on all lines were set as 100 MW. Eight interface bids were considered. For each tie-line, there were two interface bids who traded at their terminal buses but in opposite directions. The prices and maximum quantities for all interface bids were respectively set as \$0.5/(MW-h) and 100MW.

\begin{figure}[ht]
  \centering
  \includegraphics[width=0.6\textwidth]{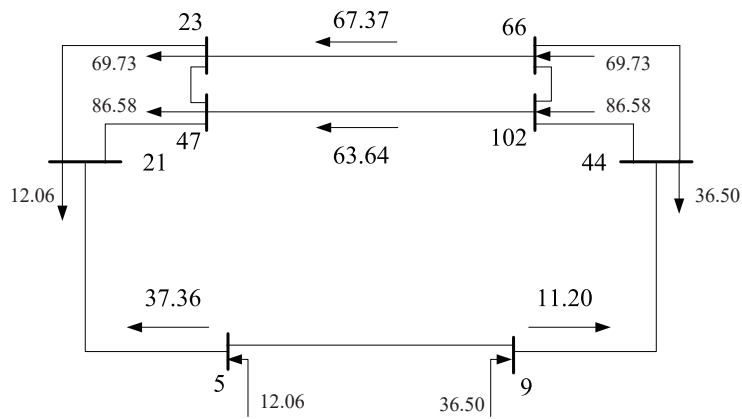}
  \caption{\small Clearing of interface bids and tie-line power flows (MW)}
  \label{fig:sys3}
\end{figure}

The results of market clearing are given in Fig. \ref{fig:sys3}, where internal parts of all areas are represented by their network equivalences. Cleared interface bids are denoted by power injections at boundary buses. Power flows through tie-lines are also shown, which were determined by the DC power flow equation for the network (\ref{eq:bnd_dclf}) in Fig. \ref{fig:sys3}.

The total cost of the three-area system was $\$1.263\times 10^{5} /h$, in which the cost of market participants was $\$601.43/h$ and the rest was the generation cost. As a reference, if there is no interchange at all, the total generation cost would be $\$1.394\times 10^{5} /h$. The reduction of generation cost largely exceeded the cost of market participants.

%

We did similar comparisons of JED, CTS, and GCTS for this three-area test as in TABLE \ref{table:compare2}. In CTS, interchange schedules were set in a pairwise manner, and proxy buses were always selected as ones with the smallest indices on their sides.

\begin{table}[ht]
\centering
\caption{\small Comparison of JED, CTS, and GCTS for the Three-area Test}
\begin{tabularx}{12cm}{cccc}
\hline
                    & JED     & CTS    & GCTS \\
\hline
\!\!Look-ahead generation costs (\$/h)  \!\! &\!\! $\!\!1.255 \!\times\! 10^5\!$  & $\!\!\!1.261 \!\times\! 10^5\!$ & $\!\!1.257 \!\times \!10^5\!$\\
\!\!Look-ahead total costs (\$/h)         & \!\!\!--      & $\!\!\!1.262 \!\times\! 10^5\!$ & $\!\!1.257 \!\times\! 10^5\!$\\
\!\!Average real-time total costs (\$/h)   \!\! & \!\!$\!1.255 \!\times\! 10^5\!$ & $\!\!\!1.263 \!\times\! 10^5\!$ & $\!\!1.263 \!\times \!10^5\!$\\
\hline
\end{tabularx}\label{table:compare2}
\end{table}

Our conclusions of comparisons were similar to those in the two-area test. GCTS had lower look-ahead costs than CTS, which was close to JED. Although its real-time costs were similar to CTS, GCTS removed the loop-flow problem in CTS. Namely, CTS suffered from overflow problems in 92 out of the 100 scenarios with randomly generated load powers.

\section{Conclusion}

The aim of this paper is to unify major approaches to interchange scheduling: JED that achieves the ultimate economic efficiency and CTS that is the state-of-the-art market solution. GCTS partially meets this goal by maintaining the same market structure as CTS while asymptotically achieving the economic efficiency of JED under given assumptions. GCTS also ensures the revenue adequacy of each system operator.

Several important issues not considered here require further investigation.  Among these are the impacts of strategic behavior of market participants, uncertainties in real-time operations, and the asynchronous mode of interchange scheduling among more than two areas.

\section{Appendix}
\subsection{Cases with More than Two Areas}
\begin{figure}[ht]
  \centering
  \includegraphics[width=0.5\textwidth]{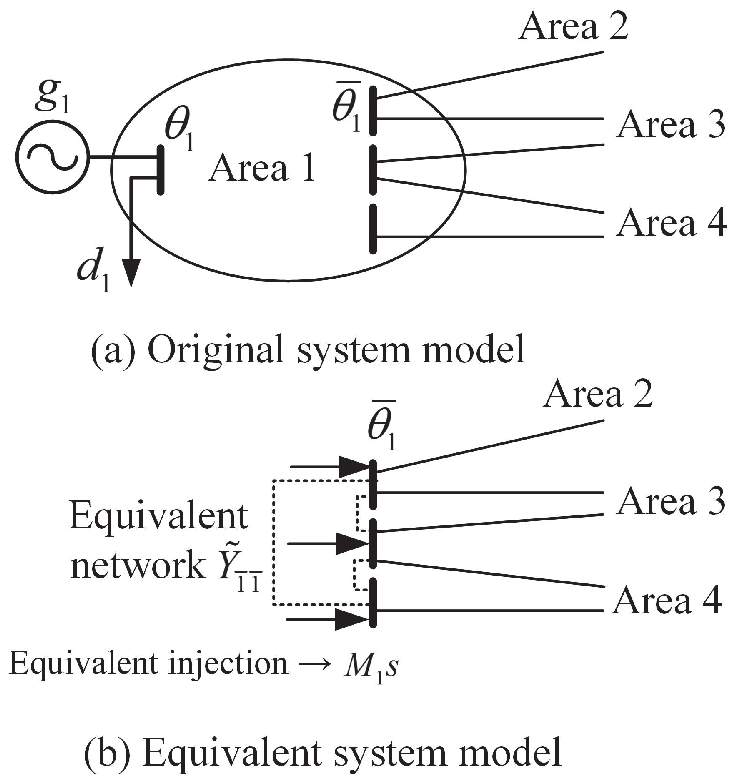}
  \caption{\small Network equivalence with more than two areas}\label{fig:threearea}
\end{figure}

In this subsection, we generalize GCTS to cases with more than two areas. For each area, the network equivalence is illustrated in Fig. \ref{fig:threearea}. Therein, internal buses are eliminated, and the equivalent admittance matrix $Y_{\bar{1}\bar{1}}$ and injection $\tilde{g}_i$ are still calculated by (\ref{eq:yeq}). The calculation of $Y_{\bar{1}\bar{1}}$ and $\tilde{g}_i$ only requires local information.

Thereby, the equivalent model of the global power system, corresponding to the Fig. \ref{fig:1}, can be obtained by eliminating all internal buses. An example of a three-area system can be seen in Fig.\ref{fig:sys3}. In the clearing of GCTS with $n$ areas, the constraint (\ref{eq:mc_bndpf}) becomes
\begin{equation}
\left[
\begin{array}{cccc}
\tilde{Y}_{\bar{1}\bar{1}}&Y_{\bar{1}\bar{2}}&\dots&Y_{\bar{1}\bar{n}}\\
Y_{\bar{2}\bar{1}}&\tilde{Y}_{\bar{2}\bar{2}}&\dots&Y_{\bar{2}\bar{n}}\\
\dots & \dots & \dots & \dots\\
Y_{\bar{n}\bar{1}}&Y_{\bar{n}\bar{2}}&\dots&\tilde{Y}_{\bar{n}\bar{n}}\\
\end{array}
\right]\!\left[\!
\begin{array}{c}
\bar{\theta}_{1}\\
\bar{\theta}_{2}\\
\dots \\
\bar{\theta}_{n}
\end{array}\!
\right]\!=\!\left[\!
\begin{array}{c}
M_1 \\
M_2 \\
\dots \\
M_n
\end{array}\!\right]s.\label{eq:threearea_bndpf}
\end{equation}

The clearing problem can be solved in a distributed fashion via existing solutions like \cite{GuoTongetc:16TPS}, which is capable to solve problems with more than two areas.

The real-time problem of Area $i$ and the settlement process are similar to those of the two-area case in Subsection III-D.

\subsection{Proof of Theorem \ref{thm:JED}}
We first show that the given assumptions ensure that the matrix $M=[M_1;M_2]$ has full row rank after removing an arbitrary phase angle reference bus. Without loss of generality, we consider bids buying from Area 1 and selling to Area 2 and assume the phase angle reference bus is in Area 1. We have assumed that there are infinitely many bids for each pair of locations when $N\rightarrow \infty$. Therefore, we can rearrange interface bids and write the matrix $M$ as follows

\begin{equation}
  \includegraphics[width=0.7\textwidth]{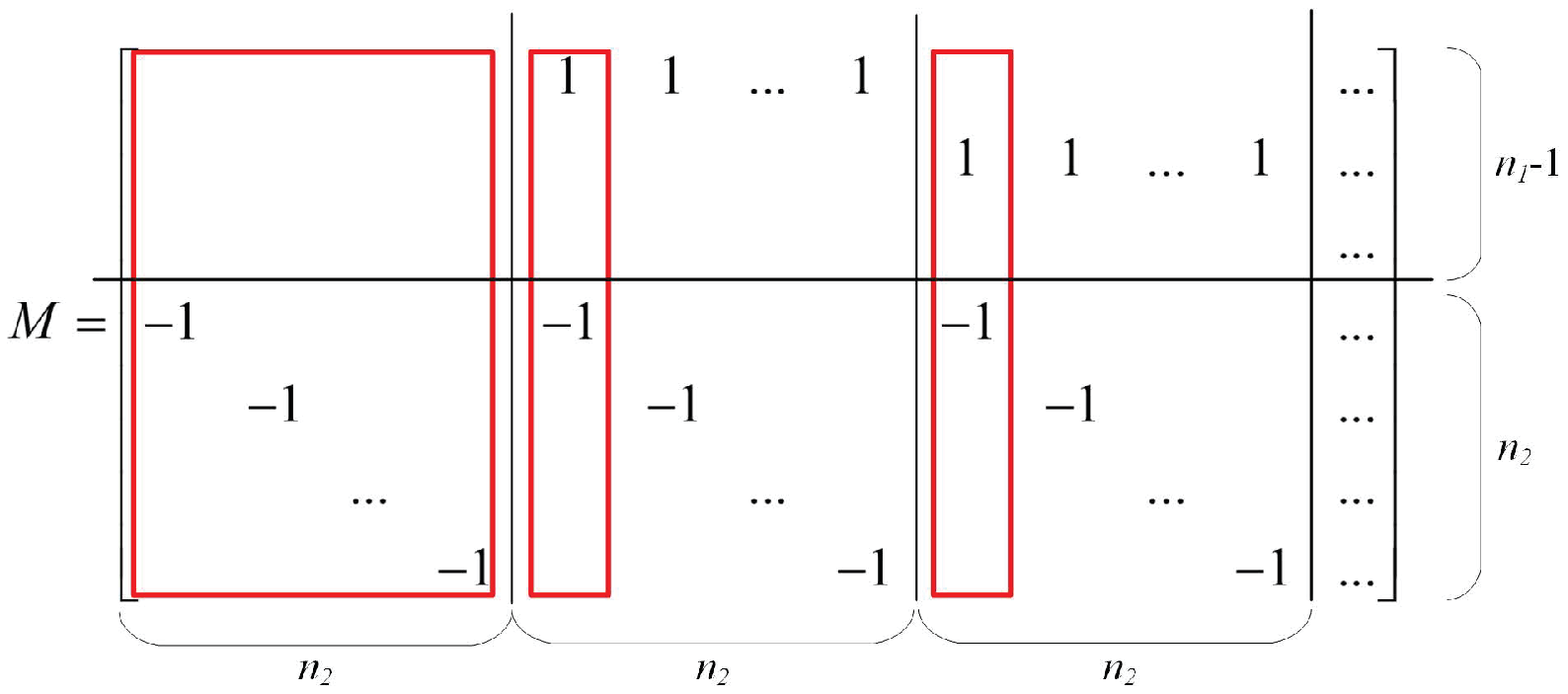}, \label{eq:Mrank}
\end{equation}
where the first block of columns includes all bids that buy at the reference bus in Area 1, the second block of columns includes all bids that buy at the first boundary bus in Area 1, etc. Scalars $n_1$ and $n_2$ are numbers of boundary buses in Area 1 and Area 2, respectively. One can pick the first block and first columns in other blocks, as highlighted by red in \eqref{eq:Mrank}, to obtain $n_1+n_2-1$ linearly independent columns. Therefore, the $M$ matrix in \eqref{eq:Mrank} has full row rank. Same conclusion holds for the other bidding direction.

Consider the optimal JED solution $\{g_{i}^{\textrm{JED}},\bar{\theta}^{\textrm{JED}}, \theta_i^{\textrm{JED}}\}$. For $N$ sufficiently large, there is an $\hat{s}$ such that the additional constraint (\ref{eq:mc_bndpf}) holds for the JED solution in GCTS. Therefore, $\{g_{i}^{\textrm{JED}},\bar{\theta}^{\textrm{JED}},$ $\theta_i^{\textrm{JED}}, \hat{s}\}$ is a feasible solution of the GCTS clearing problem (\ref{eq:mc_obj})-(\ref{eq:mc_bndpf}).

Because $\Delta \pi \rightarrow 0$ as $N\rightarrow \infty$,  we have

\begin{equation}
\lim_{N\rightarrow \infty} c(g_1^{JED},g_2^{JED},s) = \sum_i c_i(g_i^{JED}),
\end{equation}
which is the lowest possible total generation cost. The claim holds since we assume JED and GCTS are both convex programs, each having an unique global optimum.

\subsection{Proof of Theorem \ref{thm:CTS}}
When there is a single tie-line, variables $\bar{\theta}_1$ and $\bar{\theta}_2$ are scalars. The equivalent network in Fig.\ref{fig:net_eq} is simplified to Fig.\ref{fig:neteq_CTS}:

\begin{figure}[ht]
  \centering
  \includegraphics[width=0.5\textwidth]{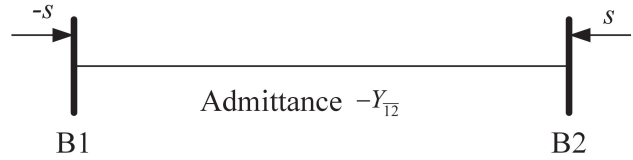}
  \caption{\small Equivalent model with single tie-line}\label{fig:neteq_CTS}
\end{figure}

Setting $\bar{\theta}_1=0$ as the phase angle reference and Area 2 exporting to Area 1 as positive, the power interchange is given by $q=-Y_{\bar{1}\bar{2}} \bar{\theta}_2$. Thereby, the GCTS clearing problem (\ref{eq:mc_obj})-(\ref{eq:mc_bndpf}) changes to
\begin{align}
\min \limits_{\{g_{i},s,\bar{\theta}_2, \theta_i, q\}}\! & c(g_1, g_2, s)=\sum \limits_{i=1}^{2}c_{i}(g_{i})+ \Delta \pi^T s, \label{eq:mc_obj_cts}\\
\textrm{subject to}&\hspace{0.1cm}H_{i}\theta_{i}+H_{\bar{i}}\bar{\theta}_{i}\leq f_{i},i=1,2, \label{eq:mc_linecons_cts}\\
&\hspace{0.1cm}q \leq \bar{f}=q_{\textrm{max}}, \label{eq:mc_tieline_cts}\\
&\hspace{0.1cm}\check{g}_i\leq g_i\leq \hat{g}_i, i=1,2,\label{eq:glimits_cts}\\
&\hspace{0.1cm} 0\leq s\leq s_{\textrm{max}}, \label{eq:slimits_cts}\\
&\hspace{-1.0cm}
\left[\!
\begin{array}{cc}
Y_{11}&Y_{1\bar{1}}\\
Y_{\bar{1}1}&Y_{\bar{1}\bar{1}}
\end{array}\!\right]\left[\!
\begin{array}{c}
\theta_{1}\\
0
\end{array}
\right]=\left[\!
\begin{array}{c}
g_{1}-d_{1}\\
q\\
\end{array}\!\right], \label{eq:mc_pf1_cts} \\
&\hspace{-1.0cm}
\left[\!
\begin{array}{cc}
Y_{22}&Y_{2\bar{2}}\\
Y_{\bar{2}2}&Y_{\bar{2}\bar{2}}
\end{array}\!
\right]\left[\!
\begin{array}{c}
\theta_{2}\\
\bar{\theta}_2\\
\end{array}
\right]=\left[\!
\begin{array}{c}
g_{2}-d_{2}\\
-q\\
\end{array}\!\right], \label{eq:mc_pf2_cts} \\
& \hspace{0.1cm} q=M_2 s.\label{eq:mc_bndpf_cts}
\end{align}
Here the matrix $M_2$ is composed of 1 and -1, depending on directions of interface bids. Note that when $q$ is fixed, the problem (\ref{eq:mc_obj_cts})-(\ref{eq:mc_bndpf_cts}) can be decoupled into three sub-problems:

i) Local economic dispatch in Area 1:
\begin{equation}
\begin{array}{ll}
\min \limits_{\{g_{1}, \theta_1\}} & c_{1}(g_{1}) \\
\textrm{subject to:} & \hspace{0.1cm}H_{1}\theta_{1}\leq f_{1},\\
&\hspace{0.1cm}\check{g}_1\leq g_1\leq \hat{g}_1,\\
&\hspace{-1.0cm}
\left[\!
\begin{array}{cc}
Y_{11}&Y_{1\bar{1}}\\
Y_{\bar{1}1}&Y_{\bar{1}\bar{1}}
\end{array}\!\right]\left[\!
\begin{array}{c}
\theta_{1}\\
0
\end{array}
\right]=\left[\!
\begin{array}{c}
g_{1}-d_{1}\\
q\\
\end{array}\!\right].\\
\end{array}\label{eq:proof_cts1}
\end{equation}

ii) Local economic dispatch in Area 2, which is similar in form with (\ref{eq:proof_cts1});

iii) The optimal clearing of interface bids given the total amount $q$:

\begin{equation}
\begin{array}{ll}
\min \limits_{s} & \Delta \pi^T s \\
\textrm{subject to:} & q=M_2 s,\\
&\hspace{0.1cm} 0\leq s\leq s_{\textrm{max}}.
\end{array}\label{eq:proof_cts2}
\end{equation}

Therefore, GCTS is to search for the optimal $q^*\in [0, q_{\textrm{max}}]$ that minimizes the sum of objective functions of the three sub-problems. Note that in Fig.\ref{fig:ctscurve}, $\pi_1 (q)$ is the derivative curve of problem (\ref{eq:proof_cts1}) with respect to $q$, $\pi_2 (q)$ is the negative derive curve of the local problem in Area 2, and $\Delta \pi (q)$ is the negative derivative of (\ref{eq:proof_cts2}). Therefore, the solution of CTS in Fig.\ref{fig:ctscurve} or Fig.\ref{fig:ctscurveCongested} is the same as that of GCTS (\ref{eq:mc_obj_cts})-(\ref{eq:mc_bndpf_cts}).

\subsection{Proof of Theorem \ref{thm:Revenue}}

For Area 1, from the optimality condition for (\ref{eq:rt_obj})-(\ref{eq:rt_lf1}) we have
\begin{equation}
\nabla_{\{\theta_1\}} L_1 = \left[
\begin{array}{cc}
Y_{11}&Y_{1\bar{1}}
\end{array}\right]\left[
\begin{array}{c}
\lambda_1^R\\
\bar{\lambda}_1^R
\end{array}\!
\right]+H_1^T \eta_1^R = 0,\label{eq:th2eq1}
\end{equation}
where $L$ is the Lagrangian of (\ref{eq:rt_obj})-(\ref{eq:rt_lf1}).

By left-multiplying $(\theta_1^*)^T$ to (\ref{eq:th2eq1}) we have
\begin{equation}
\begin{array}{ll}
(\theta_1^*)^T\nabla_{\{\theta_1\}} L_1\!&=\! (g_1\!-\!d_1)^T\lambda_1^R+f_1^T\eta_1^R - [\bar{\theta}_1^T \hspace{0.2cm}\bar{\theta}_2^T](\!\left[\!
\begin{array}{cc}
Y_{\bar{1}1}&Y_{\bar{1}\bar{1}}\\
            &Y_{\bar{2}\bar{1}}
\end{array}\!
\right]\!\!\left[\begin{array}{c}
\lambda_1^R\\
\bar{\lambda}_1^R
\end{array}\!
\right]\!+\!\left[\!
\begin{array}{c}
H_{\bar{1}}^T \eta_1^R \\
0
\end{array}\!
\right]\!)\\
&=\! (g_1\!-\!d_1)^T\lambda_1^R\!+\!f_1^T\eta_1^R\!-\!s^T M^T \!\left[\!
\begin{array}{cc}
\tilde{Y}_{\bar{1}\bar{1}}\!&\!Y_{\bar{1}\bar{2}}\\
Y_{\bar{2}\bar{1}}\!&\!\tilde{Y}_{\bar{2}\bar{2}}
\end{array}\!
\right]^{-1}\!\!\!\!\!\nabla_{\bar{\theta}} c_1^*\\
&=\! (g_1\!-\!d_1)^T\lambda_1^R\!+\!f_1^T\eta_1^R\!-s^T\mu_1=0.
\end{array}\label{eq:th2eq2}
\end{equation}
In the absence of tie-line congestion, Equation \eqref{eq:th2eq1} already proves the revenue adequacy. When tie-line congestions happen in the look-ahead dispatch, the corresponding congestion rent yield
\begin{equation}
\!\!\!
\bar{f}^T\bar{\eta}\!=\![\bar{\theta}_1^T \hspace{0.2cm}\bar{\theta}_2^T]\!\left[\!\!
\begin{array}{cc}
\tilde{Y}_{\bar{1}\bar{1}}\!\!&\!\!Y_{\bar{1}\bar{2}}\\
Y_{\bar{2}\bar{1}}\!\!&\!\!\tilde{Y}_{\bar{2}\bar{2}}
\end{array}\!\!
\right]\!\!\left[\!\!
\begin{array}{cc}
\tilde{Y}_{\bar{1}\bar{1}}\!\!&\!\!Y_{\bar{1}\bar{2}}\\
Y_{\bar{2}\bar{1}}\!\!&\!\!\tilde{Y}_{\bar{2}\bar{2}}
\end{array}\!\!
\right]^{-1}\!\!\left[\!\!\begin{array}{c}
\bar{H}_1^T\\
\bar{H}_2^T
\end{array}\!\!
\right]\!\bar{\eta}\!=\!s^T\rho.\label{eq:th2eq3}
\end{equation}
From \eqref{eq:th2eq2} and \eqref{eq:th2eq3}, we finally have
\begin{equation}
(d_1\!-\!g_1)^T\lambda_1^R\!+\!s^T(\mu_1+\frac{\rho}{2})\!=\!f_1^T\eta_1^R\!+\frac{\bar{f}^T\bar{\eta}}{2}>0. \label{eq:th2eq4}
\end{equation}

Equation \eqref{eq:th2eq4} proves Theorem \ref{thm:Revenue}. The left-hand side is the net revenue that ISO 1 collects from internal and external market participants, and the right-hand side is the sum of the internal congestion rent and a half of the tie-line congestion rent which is afforded by Area 1. Note that the net revenue in \eqref{eq:th2eq4} is non-negative because all terms on the right-hand side are non-negative.

\subsection{Proof of Theorem \ref{thm:Surplus}}
Let $\tilde{q}=1^T \tilde{s}$ and $\hat{q}=1^T \hat{s}$. We first prove that $\tilde{q}\geq\hat{q}$. Similar to the proof of Theorem \ref{thm:CTS}, the separate clearing in Area 1 is

\begin{equation}
\begin{array}{ll}
\min \limits_{\{g_{1}, \theta_1, s, q\}} & c_{1}(g_{1}) + \pi_1^T s\\
\textrm{subject to:} & \hspace{0.1cm}H_{1}\theta_{1}\leq f_{1},\\
&\hspace{0.1cm}\check{g}_1\leq g_1\leq \hat{g}_1,\\
&\hspace{-1.0cm}
\left[\!
\begin{array}{cc}
Y_{11}&Y_{1\bar{1}}\\
Y_{\bar{1}1}&Y_{\bar{1}\bar{1}}
\end{array}\!\right]\left[\!
\begin{array}{c}
\theta_{1}\\
0
\end{array}
\right]=\left[\!
\begin{array}{c}
g_{1}-d_{1}\\
q\\
\end{array}\!\right],\\
&q=M_2 s,\\
&0\leq s \leq s_{\textrm max}.
\end{array}\label{eq:proof_surplus1}
\end{equation}

Next we parameterize \eqref{eq:proof_surplus1} with respect to $q$. According to \cite{GuoBoseTong:17TPS}, the optimal value of \eqref{eq:proof_surplus1} is an piecewise affine and convex function of $q$, denoted by $J_1(q)$. We can similarly derive $J_2(q)$ for Area 2. 

Letting $q_i$ be the optimal solution of $J_i(q)$, we have $\hat{q}\leq\min\{q_1,q_2\}$. Also $\tilde{q}$ is the optimal solution of $J_1(q)+J_2(q)$. If $q_1=q_2=q$, then $\tilde{q}\geq\hat{q}$. Otherwise, we have

\begin{equation}
J_1(q_1) < J_1(\tilde{q})< J_1(q_2), \label{eq:proof_surplus2}
\end{equation}
\begin{equation}
J_2(q_2) < J_2(\tilde{q})<J_2(q_1). \label{eq:proof_surplus3}
\end{equation}

There are two possibilities that $\tilde{q}$ is smaller than both $q_1$ and $q_2$: (i) $\tilde{q}<q_2<q_1$ which violates \eqref{eq:proof_surplus2}, (ii) $\tilde{q}<q_1<q_2$ which violates \eqref{eq:proof_surplus3}. Therefore, there is always $\tilde{q}\geq\hat{q}$.

Considering the importing area whose real-time dispatch model is \eqref{eq:proof_cts1}, there is $\hat{\bar{\lambda}}_1^R >\tilde{\bar{\lambda}}_1^R$ and we have
\begin{equation}
\begin{array}{ll}
\tilde{LS}_i-\hat{LS}_i&=c_i(\hat{g}_i)-c_i(\tilde{g}_i)+\hat{\bar{\lambda}}_1^R\hat{q}-\tilde{\bar{\lambda}}_1^R\tilde{q}\\
&\geq \tilde{\bar{\lambda}}_1^R (\tilde{q}-\hat{q})+\hat{\bar{\lambda}}_1^R \hat{q}-\tilde{\bar{\lambda}}_1^R \tilde{q} \\
& = \hat{q} (\hat{\bar{\lambda}}_1^R-\tilde{\bar{\lambda}}_1^R) >0.
\end{array}
\end{equation}

Considering the exporting area whose real-time dispatch model is \eqref{eq:proof_cts1} but the term $q$ is replaced by $-q$, there is $\hat{\bar{\lambda}}_1^R < \tilde{\bar{\lambda}}_1^R$ and we have
\begin{equation}
\begin{array}{ll}
\tilde{LS}_i-\hat{LS}_i&=c_i(\hat{g}_i)-c_i(\tilde{g}_i)-\hat{\bar{\lambda}}_1^R\hat{q}+\tilde{\bar{\lambda}}_1^R\tilde{q}\\
&\geq \tilde{\bar{\lambda}}_1^R (\hat{q}-\tilde{q})-\hat{\bar{\lambda}}_1^R \hat{q}+\tilde{\bar{\lambda}}_1^R \tilde{q} \\
& = \hat{q} (\tilde{\bar{\lambda}}_1^R-\hat{\bar{\lambda}}_1^R) >0.
\end{array}
\end{equation}

There CTS achieves better local surpluses in both local markets.

{
\bibliographystyle{plain}
\bibliography{ref}
}

\end{document}